\documentclass{amsart}

\usepackage{tikz}
\usetikzlibrary{matrix,arrows}
\usepackage{geometry,latexsym,amssymb}

\theoremstyle{plain}
\newtheorem{Thm}{Theorem}
\newtheorem{Cor}{Corollary}

\newtheorem{Lem}{Lemma}

\newenvironment{Prf}{{\bf Proof:} }{\hfill $\Box$
\mbox{}}

\theoremstyle{definition}
\newtheorem{Def}{Definition}

\theoremstyle{remark}

\numberwithin{equation}{section}

\begin{document}
\title[On $G$-Sequential Continuity]%
   {On $G$-Sequential Continuity}
\author{Osman Mucuk and Tun\c{c}ar \c{S}ahan\\ Department of Mathematics, Faculty of science\\ ERC\.IYES University, Kayser\.I, Turkey}

\address{Osman MUCUK \\
           Department of Mathematics, Faculty of science,\\ Erciyes University, Kayseri, Turkey\\ +90 352 437 52 62
, +90 352 437 49 33-33208}

\email{mucuk@erciyes.edu.tr}

\address{Tun\c{c}ar \c{S}ahan\\
           Department of Mathematics, Faculty of science, Erciyes University, Kayseri, Turkey\\ +90 352 437 52 62
, +90 352 437 49 33-33212}

\email{tsahan@erciyes.edu.tr}


\date{\today}

\begin{abstract}
Let $X$ be a  first countable Hausdorff topological group.  The limit of a sequence in $X$ defines a function denoted by $lim$  from the set of all convergence sequences to $X$. This definition was modified by Connor and Grosse-Erdmann for real functions by replacing $lim$ with an arbitrary linear functional $G$ defined on a linear subspace of the vector space of all real sequences. \c{C}akall{\i} extended the concept to topological group setting and introduced  the concept of $G$-sequential compactness and investigated $G$-sequential continuity and $G$-sequential compactness in topological groups. In this paper we give a further investigation of $G$-sequential continuity in topological groups most of which are also new for the real case.
\end{abstract}

\maketitle

\section{Introduction}

The concept of continuity, and any concept related to continuity play a very important role not only in pure mathematics but also in other branches of science involving mathematics especially in computer science, information theory, biological science, and dynamical systems.

 Connor and Grosse-Erdmann \cite{ConnorGrosse} have investigated the impact of changing the definition of the convergence of sequences on the structure of sequential continuity of real functions.  \c{C}akall{\i} \cite{CakalliSequentialdefinitionsofcompactness} extended this concept to topological group setting and introduced  the concept of $G$-sequential compactness and investigated $G$-sequential continuity and $G$-sequential compactness in this generalized setting \cite{CakalliOnGcontinuity}.

 The aim of this paper is to give a further investigation of $G$-sequential continuity in topological groups and present interesting theorems which are also new for the real case.

\maketitle

\section{Preliminaries}

Before we begin some definitions and notations will be given in the following. Throughout this paper, $\textbf{N}$ will denote the set of all positive integers. Although some of the definitions that follow make sense for an arbitrary topological group, that is why we prefer using neighborhoods instead of metrics. In this paper, $X$ will always denote a topological Hausdorff group, written additively, which satisfies the first axiom of countability. We will use boldface letters $\bf{x}$, $\bf{y}$, $\bf{z}$, ... for sequences $\textbf{x}=(x_{n})$, $\textbf{y}=(y_{n})$, $\textbf{z}=(z_{n})$, ... of terms of $X$. $s(X)$, and $c(X)$ denote the set of all $X$-valued sequences, and the set of all $X$-valued convergent sequences of points in $X$, respectively.

Following the idea given in a 1946 American Mathematical Monthly
problem \cite{Buck}, a number of authors Posner \cite{Posner}, Iwinski \cite{Iwinski},
Srinivasan \cite{Srinivasan}, Antoni \cite{Antoni}, Antoni and Salat \cite{AntoniandSalat}, Spigel and
Krupnik \cite{SpielandKrupnik} have studied $A$-continuity defined by a regular
summability matrix $A$. Some authors \"{O}zt\"{u}rk \cite{Ozturk}, Sava\c{s}
and Das \cite{SavasandDas}, Sava\c{s}  \cite{Savas}, Borsik and Salat \cite{BorsikandSalat} have studied $A$-continuity for
methods of almost convergence or for related methods. See also \cite{Boos}
for an introduction to  summability matrices and  \cite{CakalliThorpe} for summability in topological groups.

The notion of statistical convergence was introduced by Fast \cite{Fast}
and has been investigated by Fridy in \cite{Fridy}. In \cite{Zymund}, Zygmund called
it "almost convergence" and established a relation between it and
strong summability. A sequence $(x_{k})$ of points in $X$ is called
to be statistically convergent to an element $\ell$ of $X$ if for each
neighborhood $U$ of $0$
\[
\lim_{n\rightarrow\infty}\frac{1}{n}|\{k\leq n: x_{k}-\ell \notin U\}|=0,
\] and this is denoted by $st-\lim_{n\rightarrow\infty}x_{n}=\ell$ (\cite{CakalliAstudyonstatisticalconvergence}). Statistical limit is an additive function on the group of statistically convergent sequences of points in $X$ (see also \cite{Cakalli2}, \cite{MaioKocinacStatisticalconvergenceintopology} and \cite{CakalliandKhan}).

A sequence $(x_{k})$ of points in a topological group is called lacunary statistically convergent to an element $\ell$ of $X$ if
\[
\lim_{r\rightarrow\infty}\frac{1}{h_{r}}|\{k\in I_{r}: x_{k}-\ell \notin U\}|=0,
\]
for every neighborhood $U$ of 0 where $I_{r}=(k_{r-1},k_{r}]$ and $k_{0}=0$, $h_{r}:k_{r}-k_{r-1}\rightarrow
\infty$ as $r\rightarrow\infty$ and $\theta=(k_{r})$ is an increasing sequence of positive integers.
For a constant lacunary sequence, $\theta=(k_{r})$, the lacunary statistically convergent sequences in a topological group form a subgroup of the group of all $X$-valued sequences, and lacunary statistical limit is an additive function on this space (see \cite{Cakalli1} for topological group setting; \cite{FridyandOrhan1}, and \cite{FridyandOrhan2} for the real case).

By a method of sequential convergence, or briefly a method, we mean an additive function $G$ defined on a subgroup $c_{G}(X)$ of $s(X)$ into $X$
\cite{CakalliSequentialdefinitionsofcompactness}. A sequence \; \; $\textbf{x}=(x_{n})$ is said to be $G$-convergent to $\ell$ if $\textbf{x}\in c_{G}(X)$ and $G(\textbf{x})=\ell$. In particular, $\lim$ denotes the limit function $\lim \textbf{x}=\lim_{n}x_{n}$ on the group $c(X)$. A method $G$ is called regular if every convergent sequence $\textbf{x}=(x_{n})$ is $G$-convergent with $G(\textbf{x})=\lim \textbf{x}$.
Clearly, if $f$ is $G$-sequentially continuous on $X$, then it is G-sequentially continuous on every subset $Z$ of $X$, but the converse is not necessarily true since in the latter case the sequences $x$ s are restricted to $Z$. This was demonstrated by an example in \cite{ConnorGrosse} for a real function.

We define sum of two methods of sequential convergence $G_{1}$ and $G_{2}$ as \[(G_{1}+G_{2})(\textbf{x})=G_{1}(\textbf{x})+G_{2}(\textbf{x})\]  where $c_{G_{1}+G_{2}}(X)=c_{G_{1}}(X)\cap {c_{G_{2}}(X)}$ (\cite{CakalliOnGcontinuity}).

The notion of regularity introduced above coincides with the classical notion of regularity for summability matrices. See \cite{Boos} for an introduction to regular summability matrices and see \cite{Zymund} for a general view of sequences of reals or complex.

 First of all, we give the definition of $G$-sequential closure of a subset of $X$. Let $A\subset X$ and $\ell \in X$. Then $\ell$ is in the $G$-sequential closure of $A$ (it is called $G$-hull of $A$ in \cite{ConnorGrosse}) if there is a sequence $\textbf{x}=(x_{n})$ of points in $A$ such that $G(\textbf{x})=\ell$. We denote $G$-sequential closure of a set $A$ by $\overline{A}^{G}$. We say that a subset $A$ is {\em $G$-sequentially closed} if it contains all of the points in its $G$-closure, i.e. a subset $A$ of $X$ is $G$-sequentially closed if $\overline{A}^{G}\subset A$.

It is clear that $\overline{\phi}^{G}=\phi$ and $\overline{X}^{G}=X$. If\; $G$ is a regular method, then $A\subset \overline{A}\subset \overline{A}^{G}$, and hence $A$ is $G$-sequentially closed if and only if $\overline{A}^{G}=A$. Even for regular methods, it is not always true that $\overline {(\overline{A}^{G})}^{G}=\overline{A}^{G}$.

Even for regular methods, the union of any two $G$-sequentially closed subsets of $X$ need not be a $G$-sequentially closed subset of $X$ as is seen by considering Counterexample 1 given after Theorem 4 in \cite{CakalliOnGcontinuity}.

\c{C}akall\i \; introduced the concept of $G$-sequential compactness and proved that $G$-sequentially continuous image of any $G$-sequentially compact subset of $X$ is also $G$-sequentially compact (see Theorem 7 in \cite{CakalliSequentialdefinitionsofcompactness}). He investigated $G$-sequential continuity, and obtained further results in \cite{CakalliOnGcontinuity}  (see also  \cite{Cakalli6}, \cite{DikandCanak}, \cite{CakalliNewkindsofcontinuities}, \cite{CakalliForwardcontinuity} and  \cite{CakalliDeltaquasiCauchysequences}  for some other types of continuities which cannot be given by any sequential method).  Among his results the following is of interest for our investigation.

\begin{Thm} \label{Thmintersecclosed}{\em (\cite{CakalliOnGcontinuity}, Theorem 5)}
 Let $G$ be a regular method and $\{A_{i}\}$ be any collection of subsets of $X$ where $I$ is any index set and $i\in{I}$. Then the following are satisfied:
\begin{description}
    \item[i]    $\bigcup_{i\in I}\overline{A_{i}}^{G} \subset \overline{\bigcup_{i\in I}A_{i}}^{G}$,
    \item[ii]    $\overline{\bigcap_{i\in I}A_{i}}^{G}\subset \bigcap_{i\in I}\overline{A_{i}}^{G}$,
    \item[iii]    $\sum_{i\in I}\overline{A_{i}}^{G}\subset \overline{\sum_{i\in I}A_{i}}^{G}$.
\end{description}
\end{Thm}

\section{Results}

In \cite{CakalliSequentialdefinitionsofcompactness} and \cite{CakalliOnGcontinuity} the concept of $G$-sequential continuity has been investigated. We give further results on $G$-sequential continuity. First we prove the following theorem.

\begin{Thm} \label{Theointersectclosed} Let $G$ be a  method on $X$. Intersection of any collections of $G$-sequentially closed subsets of $X$ is $G$-sequentially closed.
\end{Thm}

\begin{Prf}
Let $\ell\in{\overline{\bigcap_{i\in I}F_{i}}^{G}}$. Then there exists a sequence $\textbf{x}=(x_{n})$ of points in $\bigcap_{i\in I}F_{i}$ such that $G(\textbf{x})=\ell$. Hence for each $i\in I$,  $\textbf{x}$ is a sequence of points in $F_{i}$ such that $G(\textbf{x})=\ell$. This implies that $\ell\in{\overline{F_{i}}}^G$ for each $i\in{I}$. As each $F_{i}$ is $G$-sequentially closed, $\ell \in {F_{i}}$ for each $i \in{I}$. Thus $\ell \in{\bigcap_{i\in I}F_{i}}$. Hence $\overline{\bigcap_{i\in I}F_{i}}^{G} \subset{\bigcap_{i\in I}F_{i}}$. So $\bigcap_{i\in I}F_{i}$ is $G$-sequentially closed. This completes the proof of the theorem.

\end{Prf}

We again note that contrary to the expectance union of two $G$-sequentially closed subsets of $X$ need not be $G$-sequentially closed even for a regular method $G$. Define $G(\textbf{x}):= lim_{n\rightarrow\infty} (\frac{(x_{n}+x_{n+1}}{2})$, and write  $A=\{0\}$, $B=\{1\}$ where $X$ is the real space with usual topology. The sets $A=\{0\}$ and  $B=\{1\}$ are $G$-sequentially closed. But the  union of $A=\{0\}$ and  $B=\{1\}$, i.e. $A=\{0\}\cup \{1\}=\{0,1\}$ is not $G$-sequentially closed.  This example also explains that even for regular methods, it is not always true that $\overline {(\overline{A}^{G})}^{G}=\overline{A}^{G}$; since  $G$-closure of the set $A=\{0,1\}$, i.e. $\overline{A}^{G}=\{0,\frac{1}{2},1\}$, while $\overline {(\overline{A}^{G})}^{G}=\{0,\frac{1}{4},\frac{1}{2},\frac{3}{4},1 \}$ (\cite{CakalliOnGcontinuity}).

\begin{Def}
 A subset $A$ of $X$  is {\em $G$-sequentially open} if its complement is $G$-sequentially closed, i.e. $\overline{X\setminus A}^{G}\subseteq X\setminus A$.
\end{Def}

From the fact that $G$-sequential closure of a subset of $X$ includes the set itself for a regular sequential method $G$ we see that a subset $A$ is $G$-sequentially open if and only if  $\overline{X\setminus A}^{G}= X\setminus A$ for a regular sequential method $G$.

By Theorem \ref{Theointersectclosed}, we obtain the following result which states that union of any $G$-sequentially open subsets of $X$ is also $G$-sequentially open.
\begin{Thm} \label{TheoremunionGopen}  Let $G$ be a  method on $X$.  Then union of any $G$-sequentially open subsets of $X$ is $G$-sequentially open.

\end{Thm}
\begin{Prf} The proof follows easily from Theorem \ref{Theointersectclosed} so is omitted.
\end{Prf}

 Intersection of two $G$-sequentially open subsets of $X$ need not be $G$-sequentially open. Therefore  the set of $G$-sequentially open subsets of $X$ does not always give a topology on $X$.

Now we modify the ordinary concept of a neighborhood of a point to the $G$-sequential case.

\begin{Def}

Let $G$ be a  method,  $U$  a subset of $X$ and  $a\in{U}$. $U$ is called to be a {\em $G$-sequential neighborhood} of $a$ if there exists a $G$-sequential open subset $A$ of $X$ with $a\in{A}$ such that $A\subseteq{U}$.\end{Def}

It immediately follows from the definition that a subset $A$ of $X$ is $G$-sequentially open if and only if it is a neighborhood of each point of $A$.

\begin{Thm} \label{Gseqopenset} Let $G$ be a  method on $X$ and $A\subseteq X$. Then $A$ is $G$-sequentially open if and only if  each $a\in A$ has a $G$-sequentially open neighborhood $U_a$ such that $U_a\subseteq A$.
\end{Thm}
\begin{Prf} The proof is clear since  the union of $G$-sequentially open subsets is also $G$-sequentially  open by Theorem \ref{TheoremunionGopen}.
\end{Prf}
\begin{Def} \label{DefGseqinterior} Let $G$ be a  method on $X$ and $A\subseteq X$. Then the set
\[\bigcup\{U\subseteq A\mid  U~ is~ G-sequentially~open \} \]
 is called {\em $G$-interior} of $A$ and denoted by $(A^0)^G$. \end{Def}

 Here we remark that $a\in (A^0)^G$ if and only if there is a $G$-sequentially open neighborhood $U$ of $a$ such that $U\subseteq A$.

Now we give the definition of $G$-sequentially open function.

\begin{Def} A function $f$ is said to be a {\em $G$-sequentially open} if  the  image of any $G$-sequentially open subset of $X$ is $G$-sequentially open, i.e  $f(A)$ is $G$-sequentially open subset of $X$ whenever $A$ is.

\end{Def}

In \cite{CakalliOnGcontinuity} \c{C}akall{\i} gave the following definition.

\begin{Def} A function $f$ is said to be {\em $G$-sequentially closed} if the  image of any $G$-sequentially closed subset of $X$ is $G$-sequentially closed, i.e.  $f(K)$ is a $G$-sequentially closed subset of $X$ whenever $K$ is.\end{Def}
\begin{Thm} Let $G$ be a regular method on $X$.  A function $f\colon X\rightarrow X$ is $G$-sequentially closed if $\overline{f(A)}^G\subseteq f(\overline{A}^G)$ for every subset $A$.
\end{Thm}
\begin{Prf} Let $f\colon X\rightarrow X$ be a function such that $\overline{f(A)}^G\subseteq f(\overline{A}^G)$ for any  subset $A$.  Let $K$ be a $G$-closed subset. By assumption  $\overline{f(K)}^G\subseteq f(\overline{K}^G)$.  Since $G$ is regular $\overline{K}^G=K$ and so we have that $\overline{f(K)}^G\subseteq f(K)$. Therefore $f(K)$ is closed.

\end{Prf}
\begin{Thm} \label{Theointpropert} Let $G$ be a  method on $X$ and $A,B\subseteq X$. Then we have the following properties

 \begin{description}
      \item[i]   $(A^0)^G$ is $G$-sequentially open,
    \item[ii]   $(A^0)^G\subseteq A$,
    \item[iii]  $A$ is $G$-sequentially open if and only if $A=(A^0)^G$,
    \item[iv] If $A\subseteq B$, then $(A^0)^G\subseteq (B^0)^G $,
     \item[v]  $((A\cap B)^0)^G\subseteq (A^0)^G\cap (B^0)^G$,
      \item[vi] $(A^0)^G\cup (B^0)^G\subseteq ((A\cup B)^0)^G$
         \end{description}
 \end{Thm}

  \begin{Prf} (i)  $(A^0)^G$ is $G$-sequentially open as the union of $G$-sequentially open subsets included in $A$.

  (ii), (iii) and (iv)  are obvious by Definition \ref{DefGseqinterior}.

   (v) If $x\in ((A\cap B)^0)^G$, there is a $G$-sequentially open neighborhood  $U$ of $x$ such  that $x\in U\subseteq A\cap B$. This implies that $x\in U\subseteq A$ and $x\in U\subseteq B$. So $x\in (A^0)^G$ and $x\in  (B^0)^G$, i.e, $x\in (A^0)^G\cap (B^0)^G$.

 (vi) It follows from (iv).
 \end{Prf}

The $G$-sequentially interiors of arbitrary intersections and unions  are  as follows.
 \begin{Thm} \label{Leminterpropert} Let $G$ be a  method and $\{A_i\mid i\in {I}\}$  a class of subsets of $X$. Then we have the followings.

 \begin{description}
      \item[i]   $((\bigcap_{i\in I}A_i)^0)^G \subseteq \bigcap_{i\in I}(A_i ^0)^G$
       \item[ii]   $\bigcup_{i\in I}(A_i^0)^G\subseteq ((\bigcup_{i\in I} A_i)^0)^G$

               \end{description}
 \end{Thm}
 \begin{Prf} (i)  Since  $\bigcap_{i\in I}A_i \subseteq A_i$ for each $i\in I$,  we have that $((\bigcap_{i\in I}A_i)^0)^G \subseteq (A_i^0)^G$. This implies that $((\bigcap_{i\in}A_i)^0)^G \subseteq \bigcap_{i\in I}(A_i ^0)^G$.

 (ii) If $x\in{\bigcup_{i\in I}(A_i^0)^G}$, then  $x\in (A_{i_0}^0)^G$ for  an $i_0\in I$. So there exists a $G$-sequentially open neighborhood $U$ of $x$
 such that $x\in U\subseteq  A_{i_0}$. This implies that  $x\in U\subseteq \bigcup_{i\in I} A_i$  and therefore  $x\in ((\bigcup_{i\in I} A_i)^0)^G$

 \end{Prf}

 \begin{Thm}\label{Thmopenfunction}  Let $G$ be  a sequential method on  $X$.  Then a function  $f\colon X\rightarrow X$ is $G$-sequentially open if and only if
$f((A^0)^G)\subseteq ((f(A))^0)^G$  for any  subset $A\subseteq X$.
\end{Thm}
\begin{Prf}  Let the function $f\colon X\rightarrow X$ be  $G$-sequentially open  and  $A\subseteq X$.  Since $(A^{0})^G\subseteq A$, we have that   $f((A^{0})^G) \subseteq f(A)$ and therefore  $((f((A^{0})^G))^0) ^G \subseteq ((f(A))^0)^G$. Here since  $f((A^{0})^G)$  is $G$-sequentially open, it follows that    $f((A^{0})^G) \subseteq ((f(A))^0)^G$.

   Conversely suppose that  $f((A^{0})^G)\subseteq
    ((f(A))^{0})^G$ for any subset $A\subseteq X$.  So for any $G$-sequentially open subset $U$, we have that  $f(U)\subseteq (f(U)^{0})^G$ and so $f(U)$  is $G$-sequentially open. \end{Prf}

  \begin{Thm} \label{Theoremclosure} Let $G$ be a regular method and $A\subseteq X$. Then \[\overline{A}^G\subseteq \bigcap\{K\mid A\subseteq K~and ~K ~is~ G-sequentially ~closed\}.\]
 \end{Thm}
 \begin{Prf} Let $u\in \overline{A}^G$ and $K$ a $G$-sequentially closed subset including $A$. So there is a sequence $\textbf{x}=(x_n)$ of the terms in $A$ such that $G(\textbf{x})=u$. Since $\textbf{x}$  is also a sequence in $K$, we have that $u\in \overline{K}^G$. Since $K$ is $G$-sequentially closed $x\in K$. This proves the claim as required.   \end{Prf}

 \begin{Thm} \label{Theoclosure} Let $G$ be a regular method and $A\subseteq X$. If $x\in \overline{A}^G$, then for every $G$-sequentially open neighborhood $U$ of $x$, we have that $A\cap U\neq \emptyset$.
 \end{Thm}
 \begin{Prf}  Let $x\in \overline{A}^G$. Since by Theorem \ref{Theoremclosure}, $\overline{A}^G\subseteq \bigcap\{K\mid A\subseteq K~$and$ ~K ~$is$~ G-$sequentially$ ~$closed$\}$ for every $G$-sequentially closed subset containing $A$ we have  $x\in K$.  If $U$ is a $G$-sequentially open neighbourhood of $x$, then $U\cap A$ is non empty. Otherwise if $U\cap A$ is empty, then  $A\subseteq X\setminus U$ and here $X\setminus U$ is a $G$-sequentially closed subset and $x\notin X\setminus U$. This is a contradiction. \end{Prf}

We say that a subset of $X$ is {\em $G$-sequentially dense} in $X$ if $\overline{A}^G=X$.

 \begin{Cor} \label{CorGdense} If $A$ is $G$-sequentially dense in $X$, then $A\cap U$ is non empty  for each  $G$-sequentially  open subset $U$ of $X$.
 \end{Cor}

 \begin{Thm} \label{Theointpropert}  For a subset $A$ of $X$  we have that $\overline{X\setminus A}^G\subseteq X\setminus (A^0)^G$
 \end{Thm}
 \begin{Prf} If $x\in (A^0)^G$, then there is a $G$-sequentially open neighbourhood $U$ of $x$ such that $x\in U\subseteq A$.
 So $X\setminus U$ is a $G$-sequentially closed subset of $X$ containing $A$, but $x\notin  X\setminus U$. Therefore $x\notin \overline{X\setminus A}^G$ which completes the proof.
 \end{Prf}

In \cite{CakalliOnGcontinuity} the {\em boundary} of a subset $A$ is defined as the set of the points which lie in both $G$-sequential closure of $A$ and $G$-sequential closure of the complement of $A$ and denoted by $(A^{b})^{G}$.

 As a result of Theorem \ref{Theointpropert} we state   $G$-sequentially boundary of a subset  in terms of closed and open subsets as in the following.
 \begin{Cor} Let $A$ be a subset of $X$ and $(A^b)^G$   $G$-sequentially boundary of $A$. Then $(A^b)^G\subseteq \overline{A}^G\setminus (A^0)^G$.
  \end{Cor}

 \begin{Def} (\cite{ConnorGrosse} and \cite{CakalliOnGcontinuity}) A method is called {\em subsequential} if whenever $\textbf{x}$ is $G$-convergent with $G(\textbf{x})=\ell$, then there is a subsequence $(x_{n_{k}})$ of $\textbf{x}$ with $\lim_{k} x_{n_{k}}=\ell$.
\end{Def}

\begin{Def} (\cite{CakalliSequentialdefinitionsofcompactness})
A function $f\colon X \rightarrow X$ is {\em $G$-sequentially continuous} at a point $u$ if, given a sequence $(x_{n})$ of points in $X$, $G(\textbf{x})=u$ implies that $G(f(\textbf{x}))=f(u)$.
\end{Def}

\begin{Lem}\label{Lemmaclosuresame} (\cite{CakalliOnGcontinuity}) Let $G$ be a regular method. $\overline{A}^{G}=\overline{A}$ for every subset $A$ of $X$ if and only if $G$ is a subsequential method where $\overline{A}$ denotes the usual closure of the set $A$.\end{Lem}

We now prove that in the case where $G$ is a regular and subsequential method, the $G$-sequentially open subsets and usual open subsets are  same.
\begin{Lem}\label{Lemmaopensubseturesame}  If  $G$ is a regular and subsequential method, then  $(A^0)^G=A^0$ for any subset $A$ of $X$. \end{Lem}
\begin{Prf} If  $G$ is a regular and subsequential method, by Lemma \ref{Lemmaclosuresame},  $\overline{A}^{G}=\overline{A}$ for every subset $A$ of $X$.  So if $A$ is $G$-sequentially open,  $X\setminus A$ is $G$-sequentially closed and therefore by the regularity of $G$, $X\setminus A=\overline{X\setminus A}^G=\overline{X\setminus A}$. Hence $A$ is open in the usual sense. Conversely if $A$ is a usual open subset, then  $X\setminus A=\overline{X\setminus A}=\overline{X\setminus A}^G$ and so $A$ is $G$-sequentially open. Hence   $(A^0)^G$ is a usual open subset such that $(A^0)^G\subseteq A$ and so $(A^0)^G\subseteq A^0$. Similarly $A^0$ is a $G$-sequentially open subset contained by $A$ and so  $A^0\subseteq (A^0)^G$. Therefore $(A^0)^G=A^0$.
\end{Prf}

\begin{Lem} \cite{CakalliOnGcontinuity}
Let $G$ be a regular method. If a function $f$ is $G$-sequentially continuous, then $f(\overline{A}^{G})\subset{\overline{(f(A))}^{G}}$ for every subset $A$ of $X$.\end{Lem}

\begin{Lem} \cite{CakalliOnGcontinuity} Let $G$ be a regular subsequential method. Then every $G$-sequentially continuous function is continuous in the ordinary sense.\end{Lem}
Now we prove  the following Theorem.

\begin{Lem} \label{Leminverofclosed} \cite{CakalliSequentialdefinitionsofcompactness} \label{LemSeqclosed} Let $G$ be a regular method.
If a function $f$ is $G$-sequentially continuous on $X$, then inverse image $f^{-1}(K)$  of any $G$-sequentially closed subset $K$ of $X$ is $G$-sequentially closed. \end{Lem}

\begin{Thm} \label{Thminvrseimageofopen} Let $G$ be a regular method. If a function $f$ is $G$-sequentially continuous on $X$, then inverse image $f^{-1}(U)$  of any $G$-sequentially open subset $U$ of $X$ is $G$-sequentially open. \end{Thm}
\begin{Prf}
Let $f\colon X\rightarrow X$ be any $G$-sequentially continuous function and  $A$ be any $G$-sequentially open subset of $X$. Then $X\setminus A$ is $G$-sequentially closed. By Lemma \ref{Leminverofclosed},  $f^{-1}(X\setminus A)$ is $G$-sequentially closed. On the other hand
 \[f^{-1}(X\setminus A)=f^{-1}(X)\setminus f^{-1}(A)=X\setminus f^{-1}(A)\]  and so it follows that $f^{-1}(A)$ is $G$-sequentially open. This completes the proof of the theorem.

\end{Prf}

\begin{Thm} Let  $f\colon X\rightarrow X$ be a bijection. If $f$ is $G$-sequentially continuous, then
  $((f(A)) ^0)^G\subseteq f((A^0)^G)$ for any subset $A$ of $X$.
\end{Thm}
\begin{Prf} Let $A$ be a subset of $X$ and $f(A)=B$.  By  $(B^0)^G\subseteq B$, it follows that $f^{-1}((B^0)^G)\subseteq f^{-1}(B)$.  Since $f$ is a bijection $f^{-1}(B)=A$ and so we have that  $f^{-1}((B^0)^G)\subseteq A$. Then $((f^{-1}((B^0)^G))^{0})^G\subseteq (A^0)^G$. Since $f$ is sequentially continuous by Theorem \ref{Thminvrseimageofopen},  $f^{-1}((B^0)^G)$ is $G$-sequentially open and so $((f^{-1}((B^0)^G))^{0})^G=f^{-1}((B^0)^G)$.  Therefore  $f^{-1}((B^0)^G)\subseteq (A^0)^G$ and so $f(f^{-1}((B^0)^G))\subseteq f((A^0)^G)$.  Hence we obtain that  $(B^0)^G= ((f(A))^0)^G \subseteq f((A^0)^G)$.  This completes the proof of the theorem.
\end{Prf}

Now we prove that the $G$-sequential continuity of an additive function at the origin implies the $G$-sequential continuity of the function at any point in $X$, i. e. an additive function defined on $X$ to $X$ is  $G$-sequentially continuous at the origin if and only if it is   $G$-sequentially continuous at any  point $a\in X$.

\begin{Thm} \label{Thmadditivefunctioncont} Suppose that  $G$ is  a regular method. Let $f\colon X\rightarrow X$ be an additive function on $X$ into $X$. Then $f$ is $G$-sequentially continuous at the origin if and only if $f$ is $G$-sequentially continuous at any point $a\in X$.\end{Thm}

\begin{Prf} Let the additive  function $f\colon X\rightarrow X$  be $G$-sequentially continuous at the origin.  So  $G(f(\textbf{x}))=0$ whenever  $G(\textbf{x})=0$.  Let $\textbf{x}$ be a sequence in $X$ with G-lim$\textbf{x}=a$ and  $\textbf{a}$ the constant sequence $\textbf{a}=(a,a,\dots)$. Since $G$ is regular $G(\textbf{a})=a$. Therefore  the  sequence $\textbf{x}-\textbf{a}$ is $G$-convergent to the origin $0$. So by assumption $G(f(\textbf{x}-\textbf{a}))=0$.   Since $f$ and $G$ are additive $G(f(\textbf{x}))-G(f(\textbf{a}))=0$.  Here since the constant sequence $f(\textbf{a})$ tends to $f(a)$ and $G$ is regular,  $G(f(\textbf{a}))=f(a)$. Therefore we have that $G(f(\textbf{x})))=f(a)$.\end{Prf}

\begin{Cor} \label{Corgseqcont}  Let $G$ a regular method.  Then the function $f_a\colon X\rightarrow X,x\mapsto a+x$ is $G$-sequentially continuous, $G$-sequentially closed and $G$-sequentially open.\end{Cor}
\begin{Prf} Let $u\in X$ such that  $G(\textbf{x})=u$.  Then the sequence  $\textbf{a}+\textbf{x}$ is $G$-convergent to $a+u$, since the constant sequence $\textbf{a}=(a,a\dots,)$ is $G$-convergent to $a$. As $G$ is additive and regular, $G(\textbf{a}+\textbf{x})=a+u$. Therefore
\[G(f_a(\textbf{x}))=G(\textbf{a}+\textbf{x})=a+u=f_a(u)\]
and $f_a$ is $G$-sequentially continuous.

 Since the inverse of $f_a$ is $f_{-a}$, by Lemma \ref{LemSeqclosed} the function $f_a$ is $G$-sequentially closed and by   Theorem  \ref{Thmadditivefunctioncont}, $f_a$ is $G$-sequentially open.\end{Prf}

\begin{Thm} \label{TheoremGseqopentoplam} Let  $G$ be a regular sequential method. If one of the sets  $A$ and $B$ is  $G$-sequentially open, then so also is the sum $A+B$.
\end{Thm}
\begin{Prf} Suppose that $B$ is  a $G$-sequentially open subset and $A$ is any subset. By  Corollary \ref{Corgseqcont},  $a+B$ is $G$-sequentially open for any $a\in A$. Since
\[A+B=\bigcup_{a\in A}a+B\]
by Theorem \ref{TheoremunionGopen}, $A+B$ is $G$-sequentially open. \end{Prf}

\begin{Thm} \label{CorsumofGcontfunc} Let  $G$ be a sequential method on X and   $f,g\colon X\rightarrow X$ be functions on $X$. Then the following are satisfied.

 \begin{description}
    \item[i]    If $f$ and $g$ are $G$-sequentially continuous, then so also is $gf$,
    \item[ii]    If $f$ and $g$ are $G$-sequentially open (closed), then so also is $gf$,
        \item[iii]    If $f$ and $g$ are $G$-sequentially continuous, then so also is $f+g$,
     \item[iv] If $f$ and $g$ are $G$-sequentially open, then so also is $f+g$,
          \item[v] If $gf$  is    $G$-sequentially open (closed) and $f$ is onto, then  $g$ is  $G$-sequentially open (closed),
       \item[vi] If $gf$  is    $G$-sequentially open (closed) and $g$ is one to one, then  $f$ is  $G$-sequentially open (closed),

\end{description}
  \end{Thm}
\begin{Prf}
(i) Let $\textbf{x}$  be a sequence in $X$ such that $G(\textbf{x})=u\in X$.  Since  $f$ is $G$-sequentially continuous at $u$, we get $G(f(\textbf{x}))=f(u)$ and since  $g$ is $G$-sequentially continuous at $f(u)$ we have that $G(g(f(\textbf{x})))=g(f(u))$. Therefore the function  $gf$ is $G$-sequentially continuous.

(ii) is obvious.

(iii) Let $\textbf{x}$ be a sequence in $X$ with $G(\textbf{x})=u\in X$. Since the functions $f$ and $g$ are $G$-sequentially continuous,  $G(f(\textbf{x}))=f(u)$  and  $G(g(\textbf{x}))=g(u)$ . Therefore by the additivity of $G$ \[G((f+g)(\textbf{x}))=G(f(\textbf{x})+g(\textbf{x}))=G(f(\textbf{x}))+G(g(\textbf{x}))=f(u)+g(u)=(f+g)(u)\]
i.e, $f+g$ is $G$-sequentially continuous.

(iv) follows from Theorem \ref{TheoremGseqopentoplam}.

(v) Let $A$ be a $G$-sequentially open subset of $X$. Since $f$ is $G$-sequentially continuous $f^{-1}(A)$ is $G$-sequentially open. Since $gf$ is $G$-sequentially open and $f$ is onto we have that  $(gf)(f^{-1}(A))=g(A)$ is $G$-sequentially  open.  For the case where $A$ is closed, the proof is similar.

(vi) Let $A$ be a $G$-sequentially open subset of $X$. Since $gf$ is $G$-sequentially open $gf(A)$ is $G$-sequentially open.
 Since $g$ is $G$-sequentially continuous and one to one  we have that $g^{-1}gf(A)=f(A)$ is $G$-sequentially open. In the case where $A$ is closed, the proof is similar.

 \end{Prf}

The proofs of the following are straightforward but we write the details to check the conditions.
\begin{Thm} \label{TheoremcontmapGseqcont}
 Let $G$ be a sequential method. Then we have the following.
   \begin{description}
      \item[i]   If $f\colon X\rightarrow X$ is a $G$-sequentially continuous, then so  also is a restriction $f\colon A\rightarrow X$ to a subset $A$,
    \item[ii]   The identity map $f\colon X\rightarrow X$ is $G$-sequentially continuous,
    \item[iii]  For a subset $A\subseteq X$, the inclusion  map $f\colon A\rightarrow X$ is  $G$-sequentially continuous,
     \item[iv] If $G$ is regular, then the constant map  $f\colon X\rightarrow X$ is $G$-sequentially continuous,
      \item[v] If $f$  is    $G$-sequentially continuous, then so also is  $-f$,
       \item[vi] The inverse function  $f\colon X\rightarrow X,f(x)=-x$  is    $G$-sequentially continuous,
        \item[vii] The inverse function  $f\colon X\rightarrow X,f(x)=-x$  is    $G$-sequentially closed.

\end{description}
\end{Thm}

\begin{Prf} (i) Let $\textbf{x}$ be a sequence of the terms in $A$ with $G(\textbf{x})=u$. Since $f$ is $G$-sequentially continuous, $G(f(x))=f(u)$.

(ii) Let $G(\textbf{x})=u$ for a sequence $\textbf{x}$ in $X$. Then $G(f(\textbf{x}))=G(\textbf{x})=u=f(u)$ and so $f$ is $G$-sequentially continuous.

(iii) Let $\textbf{x}$ be a sequence in $A$ with $G(\textbf{x})=u\in A$. Then we have that $G(f(\textbf{x}))=G(\textbf{x})=u$.

 (iv) Let $f\colon X\rightarrow X$ be a constant map with $f(x)=x_0$ and let $\textbf{x}$ be a sequence in $X$ with $G(\textbf{x})=u$. Then $f(\textbf{x})=(x_0, x_0, ..., x_0,...)$ which converges to $x_0$. Since  $G$ is regular $G(f(\textbf{x}))=x_0=f(u)$. Therefore $f$ is $G$-sequentially continuous.

 (v) Let $\textbf{x}=(x_n)$ be a sequence in $X$ such that $G(\textbf{x})=u$.  Since $f$ is $G$-sequentially continuous $G(f(\textbf{x}))=f(u)$. Therefore  $G(-f(\textbf{x}))=-G(f(\textbf{x}))=-f(u)$ and hence $-f$ is sequentially continuous.

 (vi) For a sequence $\textbf{x}=(x_n)$ in $X$ such that $G(x)=u$, we have that $G(f(x))=G(-x)=-G(x)=-u=f(u)$. Therefore $f$ is sequentially continuous.

 (vii) Let $A$ be a $G$-sequentially closed subset of $X$. We prove that $-A$ is $G$-sequentially closed. Let $\textbf{x}$ be a sequence of the terms in $-A$ such that $G(\textbf{x})=u$.  Then  $-\textbf{x}$ is a sequence in $A$ with $G(-\textbf{x})=-u$ and therefore  $-u\in \overline{A}^G$.  Since $A$ is $G$-sequentially closed $-u\in A$.  Hence $u\in -A$ which means that $-A$ is $G$-sequentially closed.
\end{Prf}

From Theorem  \ref{CorsumofGcontfunc}  and Theorem \ref{TheoremcontmapGseqcont}, we have the following Corollary.

\begin{Cor} Let $G$ be a regular method and  $C^G(X)$  the class of $G$-sequentially continuous functions. Then $C^G(X)$ becomes a group with the sum of functions.\end{Cor}

\begin{Thm} Let  $G$ be a regular method and  $f\colon X\rightarrow X$ an additive $G$-sequentially continuous map. Then $A=\{x\in X\mid f(x)=0\}$, the kernel of $f$ is a $G$-sequentially closed subset of $X$.
\end{Thm}
\begin{Prf} If  $u\in \overline{A}^G$, there exists a sequence $\textbf{x}=(x_n)$ of the points in   $A$ such that $G(\textbf{x})=u$. So the sequence $f(\textbf{x})$ is the constant sequence with  $f(\textbf{x})=(0,0,\dots)$ and tends to $0$. Since $G$ is regular $G(f(\textbf{x}))=0$ and since $f$ is $G$-sequentially continuous, $G(f(\textbf{x}))=f(u)$. Therefore $f(u)=0$ and so $u\in A$. This proves that $A$ is $G$-sequentially closed.   \end{Prf}

\begin{Thm} Let  $G$ be a regular method and $f,g\colon X\rightarrow X$  additive $G$-sequentially continuous functions.   Then $A=\{x\in X\mid f(x)=g(x)\}$ is a $G$-sequentially closed subset of $X$.
\end{Thm}
\begin{Prf} If  $u\in \overline{A}^G$, there exists a sequence $\textbf{x}=(x_n)$ of the terms in   $A$ such that $G(\textbf{x})=u$. So the sequence $(g-f)(\textbf{x})=g(\textbf{x})-f(\textbf{x})$  is the constant sequence $(0,0,\dots)$ and tends to $0$. Since $G$ is regular, $G((g-f)(\textbf{x}))=0$ and since  $g-f$ is $G$-sequentially continuous, $G((g-f)(\textbf{x})=(g-f)(u)$. Therefore $(g-f)(u)=0$ and so $u\in A$.\end{Prf}

\section{Conclusion}

The present work improves not only the work of Connor and Grosse-Erdmann \cite{ConnorGrosse} as we have presented it in a more general setting, i.e. in a topological group which is more general than the real space, but also the papers \cite{CakalliSequentialdefinitionsofcompactness} and \cite{CakalliOnGcontinuity} of \c{C}akall{\i}, which is wholly new. So that one may expect it to be more useful tool in the field of topology in modeling various problems occurring in many areas of science, computer science, information theory, and biological science. It seems that an investigation of the present work taking "nets" instead of "sequences" could be done using the properties of "nets" instead of using the properties of "sequences". For further study, we also suggest to investigate the present work for the fuzzy case. However, due to the change in settings, the definitions and methods of proofs will not always be analogous to those of the present work (see \cite{CakalliandPratul} for the definitions in the fuzzy setting).

\end{document}